\documentclass[oneside]{amsart}

\usepackage[letterpaper,body={12.6cm,20.5cm}, mag=1000]{geometry}
\usepackage{amssymb}
\usepackage{amsthm}
\usepackage{amscd}

\numberwithin{equation}{section}
\theoremstyle{plain}
\newtheorem{cor}[equation]{Corollary}
\newtheorem{lemma}[equation]{Lemma}

\newtheorem{prop}[equation]{Proposition}

\newtheorem*{thma}{Theorem}

\theoremstyle{definition}

\newtheorem{remark}[equation]{Remark}

\newcommand{\dlabel}[1]{\ifmmode \text{\ttfamily \upshape [#1] } \else
{\ttfamily \upshape [#1] }\fi \label{#1}}\newcommand{\A}{\operatorname{A} }
\newcommand{\B}{\operatorname{B} }
\newcommand{\C}{\operatorname {C} }
\newcommand{\D}{\operatorname{D }}

\newcommand{\K}{\operatorname{K }}
\newcommand{\Ho}{\operatorname{H }}
\newcommand{\M}{\operatorname{M }}
\newcommand{\N}{\operatorname{N }}
\newcommand{\Ha}{\operatorname{H} }

\newcommand{\Z}{\operatorname{Z} }

\newcommand{\gen}[1]{\left < #1 \right >}
\newcommand{\Aut}{\operatorname{Aut} }

\newcommand{\Hom}{\operatorname{Hom} }
\newcommand{\Inn}{\operatorname{Inn} }

\newcommand{\Ker}{\operatorname{Ker} }

\newcommand{\Autcent}{\operatorname{Autcent} }

\begin{document}

\title{On central automorphisms fixing the center element-wise}

\author{Manoj K.~Yadav}

\address{School of Mathematics, Harish-Chandra Research Institute \\
Chhatnag Road, Jhunsi, Allahabad - 211 019, INDIA}

\email{myadav@hri.res.in}
\thanks{Research supported by DST (SERC Division), the Govt. of INDIA}

\date{\today}

\begin{abstract}
Let $G$ be a finite $p$-group of nilpotency class $2$. We find necessary and
sufficient conditions on $G$ such that each central automorphism of $G$ fixes
the center of $G$ element-wise.
\end{abstract}

\maketitle
\section{Introduction}

Let $G$ be a finite group. By $\Z(G)$, $\gamma_2(G)$, $\Aut(G)$ and $\Inn(G)$, 
we denote the center, the commutator subgroup,
the group of all automorphisms and the group of all inner automorphisms 
of $G$ respectively. For $x \in G$, $[G,x]$ denotes the set $\{[g,x] = 
g^{-1}x^{-1}gx | g \in G\}$.
To say that some $\Ha$ is a subgroup (proper subgroup) of $G$ we write 
$\Ha \leq G$ ($\Ha < G$). For  any group $\Ha$ and an abelian group $\K$, 
$\Hom(\Ha, \K)$ denotes the group of all homomorphisms from $\Ha$ to $\K$.
An automorphism $\alpha$ of $G$ is called \emph{central} if
$x^{-1}\alpha(x) \in \Z(G)$ for all $x \in G$. The set of all central
automorphisms of $G$, which is here denoted by $\Autcent(G)$, is a normal
subgroup of $\Aut(G)$. Notice that $\Autcent(G) = \C_{\Aut(G)}(\Inn(G))$, the
centralizer of the subgroup $\Inn(G)$ in the group $\Aut(G)$. 

Let $\M$ and $\N$ be two normal subgroups of $G$. By $\Aut^{\N}(G)$ 
 we mean the subgroup of $\Aut(G)$ consisting of all the 
automorphisms which centralize $G/\N$ and by $\Aut_{\M}(G)$ 
 we mean the subgroup of $\Aut(G)$ consisting of all the 
automorphisms which centralize $\M$. 
We denote $\Aut^{\N}(G) \cap \Aut_{\M}(G)$ by $\Aut^{\N}_{\M}(G)$. 
Throughout the paper $p$ always denotes a prime number.

In a very recent article of Attar \cite{Am07}, it is proved for a finite 
$p$-group $G$ that 
$\Aut_{\Z(G)}^{\Z(G)}(G) = \Inn(G)$ if and only if $G$ is abelian or $G$ is 
nilpotent of class $2$ and $\Z(G)$ is cyclic. So it may be interesting to
study finite $p$-groups $G$  such that $\Aut_{\Z(G)}^{\Z(G)}(G)$ contains
$\Inn(G)$ and coincides with $\Autcent(G)$. 
Notice that if $\Inn(G) \le \Autcent(G)$, then the nilpotency class of $G$ is
at most $2$. Moreover, if $G$ is abelian with 
$\Autcent(G) = \Aut_{\Z(G)}^{\Z(G)}(G)$, then
it follows that $G = 1$. So we restrict our attention to finite $p$-groups 
of class $2$. 

Let $G$ be a finite $p$-group of class $2$. 
Then $G/\Z(G)$ and $\gamma_2(G)$ have equal exponent $p^c$ (say). Let
\[G/\Z(G) = \C_{p^{a_1}} \times \C_{p^{a_2}} \times \cdots \times \C_{p^{a_r}},\]
where $\C_{p^{a_i}}$ is a cyclic group of order $p^{a_i}$, $1 \le i \le r$, and
$a_1 \ge a_2 \ge \cdots \ge a_r >0$. Let $k$ be the largest
integer between $1$ and $r$ such that $a_1 = a_2 = \cdots =a_k
= c$. Notice that $k \ge 2$. Set $\bar{\M}= \M/\Z(G) = \C_{p^{a_1}} \times  
\cdots \times \C_{p^{a_k}}$. Let 
\[G/\gamma_2(G) = \C_{p^{b_1}} \times \C_{p^{b_2}} \times \cdots \times 
\C_{p^{b_s}},\]
where $b_1 \ge b_2 \ge \cdots \ge b_s > 0$,
be a cyclic decomposition of $G/\gamma_2(G)$ such that $\bar{\M}$ is
isomorphic to a subgroup of $\bar{\N} = \N/\gamma_2(G)$ 
$:= \C_{p^{b_1}} \times \C_{p^{b_2}} \times \cdots \times \C_{p^{b_k}}$.
Using the above terminology we state our result in
\begin{thma}\label{thma}
 Let $G$ be a finite $p$-group of nilpotency class $2$. Then $\Autcent(G) = 
\Aut_{\Z(G)}^{\Z(G)}(G)$ if and only if $r = s$, $(G/\Z(G))/\bar{\M} \cong 
(G/\gamma_2(G))/\bar{\N}$ and the exponents of $\Z(G)$ and $\gamma_2(G)$ are 
equal.
\end{thma}

A generalization of the main theorem of Attar \cite{Am07} is given in
 Proposition \ref{prop1}.

As a consequence of these results, we derive the main theorem of Curran and 
McCaughan \cite{CM01} in Corollary \ref{cor1}.

\section{Proofs}

Let $G$ be a finite group and $\M$ be a central subgroup of $G$. Let $\alpha
\in \Aut^{\M}(G)$. Then we can define a homomorphism $f_{\alpha}$ from  
$G$ to $\M$ 
such that $f_{\alpha}(x) = x^{-1}\alpha(x)$. On the other hand, given a 
homomorphism
$f$ from $G$ to $\M$, we can always define an endomorphism $\alpha_f$ of $G$
such that $\alpha_f(x) = xf(x)$. But $\alpha_f$ is an automorphism of $G$ if
and only if for every non-trivial element $m \in \M$, $f(m) \neq m^{-1}$ .

The following lemma is an easy exercise.
\begin{lemma}\label{lemma0}
Let $G$ be a finite group and $\M$ be a central subgroup of $G$. Let $\M \le
\cap \Ker(f)$, where $f$ runs over all elements in $\Hom(G, \M)$.
Then the correspondence $f \rightarrow \alpha_f$ defined in the 
preceeding paragraph is a one-to-one mapping from $\Hom(G, \M)$
 onto $\Aut^{\M}(G)$. Moreover, this correspondence gives
 rise to a natural isomorphism between $\Aut_{\Z(G)}^{\M}(G)$  and 
$\Hom(G/\Z(G), \M)$.
\end{lemma}
   
A finite group $G$ is said to be \emph{purely non-abelian} if it does not 
have a non-trivial abelian direct factor. The following lemma is 
from \cite{AY65}.

\begin{lemma}\label{lemma0a}
Let $G$ be a purely non-abelian finite group. Then the correspondence $\alpha
\rightarrow f_{\alpha}$ defined above with $\M$ replaced by $\Z(G)$ is a 
one-to-one mapping of $\Autcent(G)$ onto $\Hom(G/ \gamma_2(G), \Z(G))$.
\end{lemma}

The following proposition generalizes the main result of Attar \cite{Am07}.
\begin{prop}\label{prop1}
Let $G$ be a non-abelian finite $p$-group and $\M$ be a central subgroup of 
$G$. Then $\Aut_{\Z(G)}^{\M}(G) = \Inn(G)$ if and only if the nilpotency 
class of  $G$ is $2$, $\gamma_2(G) \le \M$ and $\M$ is cyclic.
\end{prop}
\begin{proof}
Since $\M \le \Z(G)$, $f(m) = 1$ for all $f \in \Hom(G/\Z(G), \M)$ and
$m \in \M$. Thus it follows from Lemma \ref{lemma0} that there is a 
natural isomorphism  between $\Aut_{\Z(G)}^{\M}(G)$ and $\Hom(G/\Z(G), \M)$. 
Now suppose that $\Aut_{\Z(G)}^{\M}(G)$ $= \Inn(G)$. Since 
$\Inn(G) = \Aut_{\Z(G)}^{\M}(G) \le \Autcent(G)$, it follows that the 
nilpotency class of $G$ is at most $2$. As
$G$ is non-abelian, it must be of class $2$. Thus the exponents of $G/\Z(G)$
and $\gamma_2(G)$ are equal. Let $\{x_1, x_2, \cdots, x_d\}$ be a minimal
generating set for $G$. Let $i_j$ denote the inner automorphism induced by
$x_j$, i.e., $i_j(g) = x_j^{-1}g x_j$ for all $g \in G$, 
where $1 \le j \le d$. So $g^{-1}i_j(g)
\in \M$ for all $g \in G$ and $1\le j \le d$, since $\Inn(G) =
\Aut_{\Z(G)}^{\M}(G)$. But $g^{-1}i_j(g) = g^{-1}x_j^{-1}g x_j = [g,x_j]$. So 
$\{g^{-1}i_j(g) | g \in G\} = [G, x_j] \le \M$. Since 
$\gamma_2(G) = \gen{[G, x_j]| 1 \le j \le d}$, it follows that $\gamma_2(G)
\le \M$. Therefore the exponent of $\M$ is greater than or equal to the
exponent of $\gamma_2(G)$. Let the exponent of $\M$ be $p^e$. 
If possible, suppose that $\M$ is not cyclic. Then $\M = \C_{p^e} \times \N$, 
where $\C_{p^e}$ is a 
cyclic subgroup of order $p^e$ and  $\N$ is some non-trivial 
proper subgroup of $\M$. Now
\begin{align*}
|\Hom(G/\Z(G), \M)| &= |\Hom(G/\Z(G), \C_{p^e})|
|\Hom(G/\Z(G), \N)| > |G/\Z(G)|\\
& = |\Inn(G)|.
\end{align*}
Thus $|\Aut_{\Z(G)}^{\M}(G)| = |\Hom(G/\Z(G), \M)| > |\Inn(G)|.$
This contradicts the given hypothesis. Hence $\M$ must be cyclic.

Conversely suppose that the nilpotency class of $G$ is $2$, $\gamma_2(G) \le
\M$ and $\M$ is cyclic. Then 
\[|\Aut_{\Z(G)}^{\M}(G)| = |\Hom(G/\Z(G), \M)| = |G/\Z(G)| = |\Inn(G)|,\]
since the exponent of $\M$ is at least equal to the exponent of $G/\Z(G)$.
This completes the proof of the
proposition.
  \hfill $\Box$

\end{proof}

The following lemma can be proved by a counting argument using Corollary 3.3 of
\cite{CM06}. But here we prove it by a direct argument.
\begin{lemma}\label{lemma3}
Let $G$ be  a finite $p$-group such that
$\Autcent(G) = \Aut_{\Z(G)}^{\Z(G)}(G)$. Then $G$ is purely non-abelian.
\end{lemma}
\begin{proof}
Assume contrarily that $G$ is not purely non-abelian. Then $G = \Ha \times
\A$, where $\Ha$ is purely non-abelian and $\A$ is non-trivial abelian 
subgroup of $G$. 
Let $\{x_1, x_2, \cdots, x_r\}$ and $\{y_1, y_2,$ $\cdots, y_s\}$ be 
minimal generating sets for 
$\Ha$ and $\A$ respectively. Then $S := \{x_1, x_2,$ $\cdots, x_r, y_1, y_2, 
\cdots, y_s\}$ is a minimal generating set for $G$.
Obviously $\Z(\Ha) \cap \Phi(G) \neq 1$, where $\Phi(G)$ denotes the Frattini
subgroup of $G$. So we
can always choose a non-trivial element $z \in \Z(\Ha) \cap \Phi(G)$ such that 
$z^p = 1$. Define a map $f$ from $G$ to $G$ by $f(w) = wz$ for all $w \in S$.
 Now it follows from the lemma (there is only one lemma in the paper) 
of \cite{hL55} that $f$ is an automorphism of $G$.
Notice that $f \in \Autcent(G)$, but $f \not\in  \Aut_{\Z(G)}^{\Z(G)}(G)$. For,
$f(y_i) = y_iz \neq y_i$, however $y_i \in \Z(G)$ for all $i$ such that $1 \le
i \le s$. This contradicts the given hypothesis that $\Autcent(G) = 
\Aut_{\Z(G)}^{\Z(G)}(G)$. Hence $G$ must be purely non-abelian. This completes
the proof of the lemma. \hfill $\Box$

\end{proof}

\begin{lemma}\label{lemma4}
Let $\A$ and $\B$ be two finite abelian $p$-groups such that 
$\A = \C_{p^{a_1}} \times \C_{p^{a_2}}$ $\times \cdots \times \C_{p^{a_s}}$,
where $a_1 \ge a_2 \ge \cdots \ge a_s > 0$, and 
$\B = \C_{p^{b_1}} \times \C_{p^{b_2}} \times \cdots \times 
\C_{p^{b_s}}$, where $b_1 \ge b_2 \ge \cdots \ge b_s > 0$. 
Let $b_j \ge a_j$ for all $j$,
 $1 \le j \le s$, and $b_j > a_j$ for some such $j$. Let $t$ be the smallest
integer between $1$ and $s$ such that $a_j = b_j$ for all $j$ such that 
$t+1 \le j \le s$. 
Then, for any finite abelian $p$-group $\C$, $|\Hom(\A, \C)| < |\Hom(\B, \C)|$
if and only if the exponent of $\C$ is at least $p^{a_t +1}$.
\end{lemma}
\begin{proof}
Let $\Ha = \C_{p^{a_1}} \times \cdots \times \C_{p^{a_t}}$, $\K = \C_{p^{b_1}}
 \times \cdots \times \C_{p^{b_t}}$ and $\D \cong \C_{p^{a_{t+1}}} \times
 \cdots \times \C_{p^{a_s}} \cong \C_{p^{b_{t+1}}} \times
 \cdots \times \C_{p^{b_s}}$, since $b_j = a_j$ for all $j$ such that $t+1 \le
 j \le s$. Then 
\[\Hom(\A, \C) \cong \Hom(\Ha \times \D, \C) \cong \Hom(\Ha, \C) \times 
\Hom(\D, \C)\]
and 
\[\Hom(\B, \C) \cong \Hom(\K \times \D, \C) \cong \Hom(\K, \C) \times 
\Hom(\D, \C).\] 
Thus $|\Hom(\A, \C)| < |\Hom(\B, \C)|$ if and only if
$|\Hom(\Ha, \C)| < |\Hom(\K, \C)|$. Since $b_j \ge a_j$, we have 
$|\Hom(\Ha, \C)| < |\Hom(\K, \C)|$ if and only if $|\Hom(\Ha, \C_{p^c})| 
< |\Hom(\K, \C_{p^c})|$ for, at least, one cyclic group $\C_{p^c}$ 
which appears
in the cyclic decomposition of $\C$. Notice that $|\Hom(\Ha, \C_{p^c})| = 
|\Hom(\K, \C_{p^c})|$ if $c \le a_t$. Hence $|\Hom(\Ha, \C_{p^c})| 
< |\Hom(\K, \C_{p^c})|$ if and only if $c > a_t$. Since the exponent of 
$\C$ is at least $c$, this completes the proof of the lemma. \hfill $\Box$

\end{proof}

Now we are ready to prove our theorem which is stated in Section 1. In the 
proof, we use the same terminology as we used in the statement of the theorem.

\vspace{.2in}

\noindent  {\bf Proof of the theorem.} 
We have $G/\Z(G) = \C_{p^{a_1}} \times \C_{p^{a_2}} \times \cdots \times 
\C_{p^{a_r}}$ and $G/\gamma_2(G) = \C_{p^{b_1}} \times \C_{p^{b_2}} \times 
\cdots \times \C_{p^{b_s}}$. Since $G/\Z(G)$ is a quotient group of 
$G/\gamma_2(G)$, we get $r \le s$ and $b_j \ge a_j$ for all $j$
such that $1 \le j \le r$.

Suppose that $\Autcent(G) = \Aut_{\Z(G)}^{\Z(G)}(G)$. 
Then  $G$ is purely non-abelian by Lemma \ref{lemma3}. Thus by Lemma
\ref{lemma0a} there is a one-to-one correspondence between
$\Autcent(G)$ and $\Hom(G/\gamma_2(G), \Z(G))$. Also by Lemma \ref{lemma0}
(with $\M = \Z(G)$) there is a natural
isomorphism between $\Aut_{\Z(G)}^{\Z(G)}(G)$ and $\Hom(G/\Z(G), \Z(G))$. 

First we claim that $r = s$. Suppose, if possible, that $r < s$. 
Since $b_j \ge a_j > 0$ for all $j$ such that $1 \le j \le r$ and $b_j > 0$ for
$r+1 \le j \le s$, we have
\begin{align*}
|\Hom(G/\Z(G), \Z(G))| & = |\Hom(\C_{p^{a_1}} \times \cdots \times 
\C_{p^{a_r}}, \Z(G))|\\ & \le  |\Hom(\C_{p^{b_1}} \times \cdots \times 
\C_{p^{b_r}}, \Z(G))|\\
& < |\Hom(\C_{p^{b_1}} \times \cdots \times \C_{p^{b_r}}, \Z(G))|\\
& \quad \;|\Hom(\C_{p^{b_{r+1}}} \times \cdots \times \C_{p^{b_s}}, \Z(G))|\\
& = |\Hom(\C_{p^{b_1}} \times \cdots \times \C_{p^{b_s}}, \Z(G))|\\
& = |\Hom(G/\gamma_2(G), \Z(G))|.
\end{align*}
This proves that $\Aut_{\Z(G)}^{\Z(G)}(G) < \Autcent(G)$. But this contradicts
our supposition. Hence $r = s$ and our claim is true. This obviously gives the
inequalities $b_j \ge a_j$ for all $j$ such that $1 \le j \le s$. 

Now we claim that $(G/\Z(G))/\bar{\M} \cong (G/\gamma_2(G))/\bar{\N}$.
Let us suppose contrary, i.e, 
$(G/\Z(G))/\bar{\M} \not\cong (G/\gamma_2(G))/\bar{\N}$. Then $b_j > a_j$ for
some $j$ such that $k+1 \le j \le s$. For, if $b_j = a_j$ for all such $j$,
then $\A :=\C_{p^{a_{k+1}}} \times \cdots \times \C_{p^{a_s}}$
 is isomorphic to $\B := \C_{p^{b_{k+1}}} \times \cdots \times \C_{p^{b_s}}$ 
and therefore  
$(G/\Z(G))/\bar{\M} \cong \A \cong \B \cong (G/\gamma_2(G))/\bar{\N}$, which
we are not taking.
Let $t$ be the smallest integer between $k+1$ and $s$ such that $b_j = a_j$ 
for all $j$ satisfying $t+1 \le j \le s$. From our choice of $k$ we have
$p^c > p^{a_{k+1}} \ge p^{a_t}$, where $p^c$ is the exponent of both $G/\Z(G)$
as well as $\gamma_2(G)$. Since the exponent of $\Z(G)$ is greater than
or equal to $p^c$, it now follows that the 
exponent of $\Z(G)$ is at least $p^{a_t +1}$. Now applying  Lemma
\ref{lemma4} with $\C = \Z(G)$, we get 
$|\Hom(\A, \Z(G))| < |\Hom(\B, \Z(G))|$. So
\begin{align*}
|\Hom(G/\Z(G), \Z(G))| &= |\Hom(\bar{\M} \times \A, \Z(G)| = 
|\Hom(\bar{\M}, \Z(G))||\Hom(\A, \Z(G))|\\
&< |\Hom(\bar{\N}, \Z(G))||\Hom(\B, \Z(G))| 
= |\Hom(\bar{\N} \times \B, \Z(G)|\\
&= |\Hom(G/\gamma_2(G), \Z(G))|,
\end{align*}
since $|\Hom(\bar{\M}, \Z(G))| \le |\Hom(\bar{\N}, \Z(G))|$.
Thus $\Aut_{\Z(G)}^{\Z(G)}(G) < \Autcent(G)$. This again contradicts
our supposition. Hence $(G/\Z(G))/\bar{\M} \cong (G/\gamma_2(G))/\bar{\N}$.

Finally suppose that the exponent of $\Z(G)$ is not equal to the exponent of
$\gamma_2(G)$. Thus the exponent of $\Z(G)$ is at least $p^{c+1}$, since the
exponent of $\gamma_2(G)$ is $p^c$. Also $|\Z(G)| > |\gamma_2(G)|$ and 
therefore $|G/\Z(G)| < |G/\gamma_2(G)|$. 
Since $(G/\Z(G))/\bar{\M} \cong (G/\gamma_2(G))/\bar{\N}$,
it follows that $|\bar{\N}| > |\bar{\M}|$. Since $\bar{\N} = \C_{p^{b_1}}
\times \cdots \times \C_{p^{b_k}}$ and $\bar{\M} = \C_{p^{a_1}} \times 
\cdots \times \C_{p^{a_k}}$, where $b_j \ge a_j$ and $a_j = c$ for all $j$
such that $1 \le j \le k$, we conclude that $b_j > a_j$ for some such $j$.
Let $t$ be the smallest integer between $1$ and $k$ such
that $b_j = a_j = c$ for all $j$ satisfying $t+1 \le j \le k$.
We know that the exponent of $\Z(G)$ is at least $p^{c+1}= p^{a_t + 1}$.  
So applying  Lemma \ref{lemma4} with $\A = \bar{\M}$, $\B = \bar{\N}$
and $\C = \Z(G)$, we get $|\Hom(\bar{\M}, \Z(G))| < |\Hom(\bar{\N}, \Z(G))|$. 
Let $\D \cong \C_{p^{a_{k+1}}} \times \cdots \times \C_{p^{a_s}}
\cong  \C_{p^{b_{k+1}}} \times \cdots \times \C_{p^{b_s}}$.  So
\begin{align*}
|\Hom(G/\Z(G), \Z(G))| &= |\Hom(\bar{\M} \times \D, \Z(G))| = 
|\Hom(\bar{\M}, \Z(G))||\Hom(\D, \Z(G))|\\
&< |\Hom(\bar{\N}, \Z(G))||\Hom(\D, \Z(G))| 
= |\Hom(\bar{\N} \times \D, \Z(G))|\\
&= |\Hom(G/\gamma_2(G), \Z(G))|.
\end{align*}
Thus $\Aut_{\Z(G)}^{\Z(G)}(G) < \Autcent(G)$, which is again a contradiction.
Hence the exponents of $\Z(G)$ and $\gamma_2(G)$ are
equal. This completes the necessary part of the theorem.

Conversely suppose that $r = s$, $(G/\Z(G))/\bar{\M} \cong 
(G/\gamma_2(G))/\bar{\N}$ and the exponents of $\Z(G)$ and $\gamma_2(G)$ 
are equal.
Since $r = s$, it follows that $G$ is purely non-abelian. For, suppose that
$G$ is not purely non-abelian. Then $G = \Ho \times \A$, where $\A$ is
non-trivial abelian and $\Ho$ is purely non-abelian subgroup of $G$. 
So $\gamma_2(G) = \gamma_2(\Ho)$ and
$\Z(G) = \Z(\Ho) \times \A$. This implies that 
\begin{eqnarray*}
G/\gamma_2(G) &=& G/\gamma_2(\Ho) = (\Ho \times \A)/\gamma_2(\Ho) 
\cong (\Ho/\gamma_2(\Ho)) \times (\A \gamma_2(\Ho)/\gamma_2(\Ho))\\
&\cong& (\Ho/\gamma_2(\Ho)) \times \A
\end{eqnarray*}
and
\begin{eqnarray*}
G/\Z(G) &=&  (\Ho \times \A)/(\Z(\Ho) \times \A) = \Ho \Z(\Ho)A/\Z(\Ho)A \cong
\Ho/(\Ho \cap \Z(\Ho)A)\\
&=& \Ho/\Z(\Ho).
\end{eqnarray*}
These equations imply that the rank of $G/\gamma_2(G)$ is strictly greater
than the rank of $\Ho/\gamma_2(\Ho)$ and  the rank of $G/\Z(G)$ is equal to the
rank of $\Ho/\Z(\Ho)$, where the rank of a finite abelian $p$-group $X$ is
defined to be the number of non-trivial cyclic factors in a cyclic
decomposition of $X$.
Since the nilpotency class of $\Ho$ is $2$, we have $\gamma_2(\Ho) \subseteq
\Z(\Ho)$. So $\Ho/\Z(\Ho)$ is a quotient group of  $\Ho/\gamma_2(\Ho)$. Thus
it follows that the rank of $\Ho/\Z(\Ho)$ is less than or equal to the rank of 
$\Ho/\gamma_2(\Ho)$. Hence the rank of $G/\gamma_2(G)$ is strictly greater
than the rank of $G/\Z(G)$. This implies that $s > r$, because $s$ is the
rank of $G/\gamma_2(G)$ and $r$ is the rank of $G/\Z(G)$. This contradition
proves that $G$ is purely non-abelian.

Since $r = s$ and $(G/\Z(G))/\bar{\M} \cong (G/\gamma_2(G))/\bar{\N}$, we have
$G/\Z(G) = \C_{p^{a_1}} \times \C_{p^{a_2}}$ 
$\times \cdots \times \C_{p^{a_s}}$
and $G/\gamma_2(G) = \C_{p^{b_1}} \times \C_{p^{b_2}} \times \cdots \times 
\C_{p^{b_s}}$ such that $\C_{p^{a_{k+1}}} \times  \cdots \times \C_{p^{a_s}}$
$\cong \C_{p^{b_{k+1}}} \times  \cdots \times \C_{p^{b_s}} \cong \D$ (say). So
the smallest $t$, $1 \le t \le s$ such that $a_j = b_j$ for all $j$ satisfying
 $t+1 \le j \le s$, is at the most $k$. 

Now suppose that $\Aut_{\Z(G)}^{\Z(G)}(G) < \Autcent(G)$. 
Then by Lemma \ref{lemma0} (with $\M = \Z(G)$) and Lemma \ref{lemma0a}, 
we have
\[|\Hom(G/\Z(G), \Z(G))| < |\Hom(G/\gamma_2(G), \Z(G))|,\]
since $G$ is purely non-abelian.
So it follows that $b_j > a_j$ for some $j$ such that $1 \le j \le s$ (to be
more precise, it happens for some $j$ such that $1 \le j \le k$). For, if $b_j
= a_j$ for all $j$, then 
$|\Hom(G/\Z(G), \Z(G))| = |\Hom(G/\gamma_2(G), \Z(G))|$, which we are not 
taking. Now applying
Lemma \ref{lemma4} with $\A = G/\Z(G)$, $\B = G/\gamma_2(G)$ and $\C = \Z(G)$,
we conclude that the exponent of $\Z(G)$ must be at least $p^{a_t +1}$,
 which is
strictly bigger than $p^{a_k} = p^c$ - the exponent of $\gamma_2(G)$. This
contradicts our supposition that $\Z(G)$ and $\gamma_2(G)$ have equal
exponents. Hence $\Aut_{\Z(G)}^{\Z(G)}(G) = \Autcent(G)$, 
which completes the proof of the  theorem.
\hfill $\Box$

\vspace{.2in}

It is fairly easy to deduce the following corollary from Proposition
\ref{prop1} and the theorem. 
\begin{cor}\label{cor1}
Let $G$ be a non-abelian finite $p$-group. Then $\Autcent(G) = \Inn(G)$ if and
only if $\Z(G) = \gamma_2(G)$ and $\Z(G)$ is cyclic.
\end{cor}

\begin{remark}
It is of interest to find  necessary and sufficient conditions for a
finite $p$-group $G$ of arbitrary nilpotency class such that $\Autcent(G) = 
\Aut_{\Z(G)}^{\Z(G)}(G)$.
\end{remark}

\noindent{\bf Acknowledgements.} I thank Prof. E. C. Dade for some useful
discussion and the referee for his/her useful comments and suggestions.

\end{document}